\documentstyle[12pt]{article}
\catcode`\@=11
\@addtoreset{equation}{section}

\catcode`\@=12
\newtheorem{Theorem}{Theorem}[section]

\newtheorem{Corollary}{Corollary}[section]

\title{A Proof On Arnold Chord Conjecture}

\author{Renyi Ma \\
Department of Mathematical Sciences \\
Tsinghua University \\
Beijing, 100084\\
People's Republic of China\\
rma@math.tsinghua.edu.cn}

\date { }

\begin{document}
\textwidth=165mm
\textheight=185mm
\parindent=8mm
\frenchspacing
\maketitle

\begin{abstract}
In this article, we first give a proof on the Arnold chord conjecture
which states that every Reeb flow has at least as many Reeb chords
as a smooth function on the Legendre submanifold has critical
points on contact manifold. Second, we prove that every Reeb flow has at least as
many close Reeb orbits as a smooth round function on the close
contact manifold has critical circles on contact manifold.
This also implies a proof on
the fact that there exists at least number $n$ close Reeb orbits on
close $(2n-1)$-dimensional convex hypersurface in $R^{2n}$ conjectured by
Ekeland.
\end{abstract}
\noindent{\bf Keywords} Contact manifold, Reeb Chord, Periodic
orbits.

\noindent{\bf 2000 MR Subject Classification} 32Q65, 53D35,53D12

\section{Introduction and results}

Let $\Sigma$ be a smooth closed oriented manifold of dimension
$2n-1$. A contact form on $\Sigma$ is a $1-$form such that $\lambda
\wedge (d\lambda )^{n-1}$ is a volume form on $\Sigma$. Associated
to $\lambda$ there are two important structures. First of all the
so-called Reeb vectorfield $\dot x=X$ defined by
$$i_X\lambda   \equiv 1, \ \ i_Xd\lambda  \equiv 0;$$
and secondly the contact structure $\xi =\xi _{\lambda }
\mapsto \Sigma $ given by
$$\xi _{\lambda }=\ker (\lambda )\subset T \Sigma .$$
By a result of Gray, \cite{gra}, the contact structure is very
stable. In fact, if $(\lambda  _t  )_{t\in  [0,1]}$ is a smooth
arc of contact forms inducing the arc of contact structures
$(\xi _t)_{t\in [0,1]}$, there exists a smooth arc
$(\psi _t)_{t\in [0,1]}$
of diffeomorphisms with $\psi _0=Id$, such that
\begin{equation}
T\Psi _t(\xi _0)=\xi _t  \label{eq:1.1}
\end{equation}
here it is important that $\Sigma $ is compact. From (\ref{eq:1.1})
and the fact that $\Psi _0=Id$ it follows immediately that there
exists a smooth family of maps $[0,1]\times \Sigma \mapsto (0,\infty
):(t, m)\to f_t(m)$ such that
\begin{equation}
\Psi ^*_t\lambda _t=f_t\lambda _0
\end{equation}
In contrast to the contact structure the dynamics of the Reeb
vectorfield changes drastically under small perturbation and in
general the flows associated to $X_t$ and $X_s$ for $t\neq s$ will
not be conjugated.

Let $(\Sigma ,\lambda )$ is a contact manifold with contact form
$\lambda $ of dimension $2n-1$, then a Legendre submanifold is a
submanifold ${\cal L}$ of $\Sigma $, which is $(n-1)$-dimensional
and everywhere tangent to the contact structure $\ker \lambda $.
Then a characteristic chord for $(\lambda ,{{\cal {L}}})$ is a
smooth path
$$x:[0,T]\to M,T>0$$
with
$$\dot x(t)=X_{\lambda }(x(t)) \ for \ t\in(0,T),$$
$$x(0),x(T)\in {\cal {L}}$$

  The main results of this paper is following:

\begin{Theorem}
Let $(\Sigma ,\lambda )$ be a contact manifold
with contact form
$\lambda $, $X_{\lambda } $ its Reeb vector field, ${\cal {L}}$ a
closed Legendre submanifold. Then there exists at least as many Reeb
characteristic chords for $(X_\lambda ,{\cal{L}})$ as a smooth
function on the Legendre submanifold ${\cal {L}}$ has critical
points.
\end{Theorem}
Theorem1.1 is asked in \cite{ar}. In \cite{ma1}, we have proved there exists at least one Reeb chord
by Gromov's $J-$holomorphic curves. Partial results is obtained in
\cite{ma,mo}.

   We recall that a round function $f$ on $M$ is the one whose critical sets
consist of smooth circles $\{C_i,i=1,...,k\}$
(see\cite{as}). One can define Round Morse-Bott function.

\begin{Theorem}
Let $(\Sigma ,\lambda )$ be a contact manifold
with contact form
$\lambda $, $X_{\lambda } $ its Reeb vector field.
Then there
exists at least as many close Reeb characteristic orbits for
$X_\lambda $ as a smooth round function on the $\Sigma $ has
critical circles.
\end{Theorem}
In \cite{ma2}, we have proved there exists at least one close
Reeb orbit by Gromov's $J-$holomorphic curves.
Theorem1.2 is related to the Arnold-Ginzburg question on magnetic field and
partial results was obtained in \cite{ar,ba,gi}.

\begin{Corollary}
Let $(\Sigma ,\lambda )$ be a close $(2n-1)-$dimensional star-shaped
hypersurface in $R^{2n}$ with contact form $\lambda $, $X_{\lambda } $ its Reeb
vector field. Then there exists at least number $n$ close Reeb
orbits.
\end{Corollary}
This implies that the Ekeland conjecture holds in \cite{ek}.

The proofs of Theorem1.1-1.2 are the extension of the methods in
\cite{ma3}.

\section{Proof of Theorem1.1}

{\bf Proof of Theorem1.1:} Let $(\Sigma ,\lambda )$ be a contact manifold
with the contact form
$\lambda $. By Whitney's embedding theorem, we first embed $\Sigma $ in $R^N$, then by considering  the cotangent bundles, the symplectizations and contactizations, one can embed $(\Sigma ,\lambda )$ into $(S^{2N+1},\lambda _0)$ with $\lambda =f\lambda _0$,
$f$ is positive function on $\Sigma $.
By the contact tubular neighbourhood theorem, the neighbourhood $U(\Sigma )$ is
contactomorphic to the symplectic vector bundle $E$ on $\Sigma $ with symplectic fibre
$(R^{2N-2n+2},\omega _0)$. By our construction, it is easy to see that there exists
Lagrangian sub-bundle $L$ in $E$ with Lagrangian fibre $R^{N-n+1}$. So, we can
extend the contact form $\lambda $ on $\Sigma $ to the neighbourhood $U(\Sigma )$
as $\bar \lambda $ such that the Reeb vector fields $X_{\bar \lambda }|\Sigma =X_\lambda $.

  We extend $f$ positively to whole
$S^{2N+1}$. So the contact form $\bar \lambda $ on $U(\Sigma )$ is extended to whole
$S^{2N+1}$ as $f\lambda _0$.

  Let $\lambda _s=(sf+1-s)\lambda _0$ be the one parameter family of contact forms
on $S^{2N+1}$. Let $X_{\lambda _s} $
its Reeb vector field.

    Let ${\cal L}$ be a close Legendre submanifold contained
in $\Sigma $, i.e., $T{{\cal L}}\subset \xi$.
Let $\bar {{\cal L}}$ be a Legendre submanifold contained
in $U(\Sigma )$ which is fibred on ${\cal L}$ and $T\bar {{\cal L}}\subset \xi _0$,
here $\xi _0$ is the standard contact structure on $(S^{2N+1},\lambda _0)$.

    Let $S$ be a smooth
hypersurface in $\Sigma $ which contains ${\cal L}$. We can assume that
the Reeb vector fields $X_{\lambda } $ is transversal to $S$.
Let $\bar S$ be a smooth
hypersurface in 
$S^{2N+1}$
which contains $\bar {\cal L}$ and extends $S$. We can assume that
the Reeb vector fields $X_{\lambda _s} $ is transversal to $\bar S$.

    Let $\eta _{s,t}$ be the Reeb flow generated by the Reeb vector field
$X_{\lambda _s}$. Let $f_s:\bar {\cal L}_s\to R$ be the arrival time function of
the Reeb flow
$\eta _{s,t}$ from the parts $\bar {\cal {L}}_s$ of $\bar {\cal {L}}$ to $\bar S$.
Then $\varphi _s=\eta _{s,f_s}(\cdot ):\bar {\cal L}_s\to \bar S$ is an exact Lagrange
embedding for symplectic form $d\lambda _s|\bar S$.

   One observes that $\varphi _s(\bar {{\cal {L}}}_s)\cap \bar {\cal {L}}$ corresponds to the Reeb chords
of $X_{\lambda _s}$.

    Now we extend the family of functions $f_s$ on $\bar {{\cal {L}}}_s$ to the
whole $\bar {{\cal {L}}}$.
This extends
$\varphi _s=\eta _{s,f_s}(\cdot ):\bar {\cal L}\to S^{2N+1}$ as an exact isotropic
embedding for symplectic form $d\lambda _s$ on $S^{2N+1}$. Then, we obtain
${\bar {\varphi }}_s=\eta _{s,f_s+t}(\cdot ):\bar {\cal L}\times [-\varepsilon ,\varepsilon ]\to S^{2N+1}$ as an exact Lagrangian
embedding for symplectic form $d(e^t\lambda _s)$ on $R\times S^{2N+1}$.
By Moser's stability theorem, this defines
an exact Lagrangian isotopy ${\bar {\varphi }}_s:\bar {\cal L}\times [-\varepsilon ,\varepsilon ]\to (R\times S^{2N+1}, d(e^t\lambda _0)$. It is well known that
an exact Lagrangian isotopy is determined by Hamilton isotopy with hamilton function
$h_s$ on $R\times S^{2N+1}$. Let $X_{h_s}$ be the hamilton vector field induced by the
hamilton function $h_s$.

   We first assume that near the Lagrangian submanifold
${\cal L}$ in $S$ the form
$\lambda $ is the standard Liouville form in
$T^*{\cal L}$. Also, We assume that near the Lagrangian submanifold
$\bar {\cal L}\times [-\varepsilon ,\varepsilon ]$ the form
$e^t\lambda _0$ is the standard Liouville form in
$T^*(\bar {\cal L}\times [-\varepsilon ,\varepsilon ])$.
Consider the Lie derivative  $L_{X_{h_s}}e^t\lambda _0$,
by hamilton perturbation, we can assume that
$L_{X_{h_s}}e^t\lambda _0=e^t\lambda _0$
on ${\bar {\varphi }}_s({\cal L}\times [-\varepsilon ,\varepsilon ])$.
Moreover, $L_{X_{h_s}}e^t\lambda _0=dH_s$, here 
$H_s$ is defined on $R\times S^{2N+1}$. The level sets of 
$H_s$ defines a foliation ${\cal {F}}_s$ on 
$R\times S^{2N+1}$. For $s$ is small enough, the Lagrangian submanifold 
$\varphi _s({\cal L})$ in $S$ is transversal 
to the foliation ${\cal {F}}_s|S$ except the Lagrangian submanifold 
${\cal L}$. So, we can perturb the 
Lagrangian isotopy ${\bar {\varphi }}_s$ such that the intersection points sets 
$\bar \varphi _s({\cal L})\cap {\cal L}$ is invariant and
the isotropic submanifold
${\bar {\varphi }}_s({\cal L})$ is transversal
to the foliation ${\cal {F}}_s$ except the intersection points set 
$\bar \varphi _s({\cal L})\cap {\cal L}$.

    This shows that the critical points of $H_1|\varphi _1({\cal L})$ are 
in the intersection points set
$\varphi _1({\cal L})\cap {\cal L}$.

\vskip 3pt

This yields Theorem1.1.

\section{Proof of Theorem1.2}

{\bf Proof of Theorem1.2:} Let $(\Sigma ,\lambda )$ be a close
contact manifold with contact form $\lambda $ and $X_\lambda $ its
Reeb vector field, then $X_\lambda $ integrates to a Reeb flow $\eta
_t$ for $t\in R^1$. Consider the form $d(-e^a\lambda )$ at the point
$(a,\sigma )$ on the manifold $(R\times \Sigma )$, then one can
check that $d(-e^a\lambda )$ is a symplectic form on $R\times \Sigma
$. Let
$$(\Sigma ',\lambda ')=((R\times \Sigma )\times \Sigma ,
(-e^a\lambda )\oplus \lambda ).$$ Then $(\Sigma ',\lambda ')$ is a
non-compact contact manifold. But the Reeb flow $\eta '_s=(id,\eta
_s)$ of Reeb vector field $X_{\lambda '}$ behaves as close contact
manifold, especially the compact set $([-a,a]\times \Sigma )\times
\Sigma $ is invariant under Reeb flow. Let
$${\cal {L}}=\{
((0,\sigma ),\sigma )|\sigma \in \Sigma \}.$$ Then, ${\cal {L}}$ is
a close Legendre submanifold in $(\Sigma ',\lambda ')$.

Then, by considering the round function and round intersection, the
method of proving Theorem1.1 yields Theorem1.2. e.q.d.

\end{document}